\date{}
\title {One-arm exponent for critical 2D percolation}
\author {Gregory F.~Lawler\thanks{Duke University and Cornell
University, Research
supported by the National Science Foundation}
\and Oded Schramm\thanks{Microsoft Research}
\and Wendelin Werner\thanks{Universit\'e Paris-Sud}
}
\documentclass[12pt,naturalnames]{article}
\usepackage{amsmath}
\usepackage{amsthm}
\usepackage{amsfonts}
\usepackage{graphicx}
\input labelfig.tex
\newif\ifhyper\IfFileExists{hyperref.sty}{\hypertrue}{\hyperfalse}
\ifhyper\usepackage{hyperref}\fi

\newif\ifdraft
\drafttrue
\long\def\note#1/{{\bf[#1] }}

\newtheorem{theorem}{Theorem}
\numberwithin{theorem}{section}
\numberwithin{equation}{section}

\newtheorem{lemma}[theorem]{Lemma}

\def\eref#1{(\ref{#1})}
\def\Prob{{\bf P}}
\def \E {{\bf E}}
\newcommand{\R}{\mathbb{R}}
\newcommand{\C}{\mathbb{C}}

\newcommand {\CS}{{C}}
\def\kk{G}
\def\H{\mathbb{H}}
\newcommand{\U}{\mathbb{U}}
\def\md{\mid}
\def\Bb#1#2{{\def\md{\bigm| }#1\bigl[#2\bigr]}}
\def\BB#1#2{{\def\md{\Bigm| }#1\Bigl[#2\Bigr]}}

\def\Pb{\Bb\P}
\def\Eb{\Bb\E}
\def\PB{\BB\P}

\def \P {{\bf P}}
\def \p {{\partial}}
\def\closure{\overline}
\def\dist{\mathrm{dist}}

\def \proof {{ \medbreak \noindent {\bf Proof.} }}
\def\proofof#1{{ \medbreak \noindent {\bf Proof of #1.} }}

\def\QED{\qed\medskip}
\def\sle#1/{$\mathrm{SLE}_{#1}$}
\def\ch{{\mathrm{ch}}}
\def\ra{{\mathrm{ra}}}
\def\rh{{\mathrm{rhs}}}
\def\th{{\tilde h}}
\def\Re{\mathop{\mathrm{Re}}}

\def\diam{\mathop{\mathrm{diam}}}
\def\crad{\mathfrak r}

\begin {document}

\maketitle

\begin {abstract}
The probability that the 
cluster of the origin in critical site percolation
on the triangular grid has diameter larger than $R$
is proved to decay like $R^{-5/48}$ as $R\to\infty$.
\end {abstract}

\section {Introduction}
Critical site percolation on the triangular lattice
is obtained by declaring independently 
each site to be open
with probability $p=1/2$ and otherwise to be closed.
Let $0\leftrightarrow\CS_R$ denote the event that
there exists an open path from
the origin to the circle $\CS_R$ of radius $R$ around the
origin.
In this paper we prove:
\begin{theorem}\label{t.pexpo}
For critical site percolation on
the standard triangular grid in the plane,
$$
\Pb{ 0 \leftrightarrow \CS_R } = R^{-5/48+ o(1)},\qquad R\to\infty\,.
$$
\end{theorem}
This result is based on the recent proof by Stanislav
Smirnov~\cite{Spaper}
that the scaling limit of this percolation process exists 
and is described by \sle6/, the stochastic Loewner evolution with
parameter $\kappa=6$.
This latter result has been conjectured in \cite {S} where
\sle\kappa/ was introduced.
See~\cite{SmW} for a treatment of other critical exponents for 
site percolation on the triangular lattice.

The exponent determined in Theorem \ref {t.pexpo}
is sometimes called the one-arm exponent.
In Appendix~\ref{a.twoarm} we briefly discuss what the methods of
the present paper can say about the  monochromatic two-arm
or backbone exponent,  
for which no established conjecture existed.
More precisely, we 
show that this exponent is the highest eigenvalue of 
a certain differential operator with mixed boundary
conditions in a triangle. 

Theorem~\ref{t.pexpo} has been conjectured in the
theoretical physics literature \cite {DN,N,N2,ADA},
 usually
in forms involving other 
critical exponents that imply 
this one using the scaling relations that have been proved
mathematically by
Kesten \cite {K}.  It is conjectured
that the theorem holds for any planar lattice.  However,
Smirnov's results mentioned above have
been established only for site percolation on the triangular lattice,
and hence we can only prove our result in this case.
In the literature, the exponent calculated in Theorem~\ref{t.pexpo}
is sometimes denoted by $1/\rho$.

Theorem~\ref{t.pexpo} will appear as a corollary 
of the following similar
result about the exponent for the scaling limit.
The scaling limit is the process obtained from percolation
by letting the mesh of the grid tend to zero
(see the next section for more detail).

\begin{theorem}\label{t.contexpo}
Let $Q$ denote the union of the clusters meeting the
unit circle $\p\U$ in the scaling limit of critical site
percolation on  the triangular lattice.  There is
a constant $c>0$ such that for all $r\in(0,1/2)$,
$$
 c^{-1} \,r^{5/48}\le
\Pb{\dist(Q,0)<r}\le c \,r^{5/48}.
$$
\end{theorem}

This theorem resembles in statement and in proof the
determination of the Brownian motion exponents carried
out in~\cite{LSW2}.
Also related to this theorem is the
calculation in~\cite{Spercform} of the probability
of an event for the scaling limit.

\medskip

We will assume that the reader is familiar with  the theory of
percolation in the plane (mainly, the Russo-Seymour-Welsh 
theorem and its consequences),
such as appearing in~\cite{G,Kbook}.  Also,
familiarity with SLE will be assumed.
To learn about the basics of SLE
the reader is advised to consult the first few sections
in~\cite{LSW1,LSW2,RS,L}.

\medskip

We now turn to discuss the SLE counterpart of Theorem~\ref{t.contexpo}.
Let
\begin{equation}  \label{e.lqdef}
 \lambda = \lambda(\kappa) := \frac{\kappa^2 - 16}{32 \kappa}\,.
\end{equation}

\begin {theorem} \label{t.sletwist}
For every $\kappa >4$ with $\kappa \neq 8$,
there  exists a constant $c>0$ such the radial \sle\kappa/ path
$\gamma:[0,\infty)\to\closure\U$ satisfies for all $r \in (0,1)$,
$$ c^{-1}
r^{\lambda} \le
\Pb{ \gamma [0,T_r]
\hbox { contains no counterclockwise loop around } 0 }
\le c\, r^{ \lambda } 
$$
where $T_r$ denotes the first time at which $\gamma$ intersects
the circle of radius $r$ around the origin.
\end {theorem}

The existence of a counterclockwise 
loop around $0$ in $\gamma[0,t]$ means
that there are $0\le t_1\le t_2\le t$ such that
$\gamma(t_1)=\gamma(t_2)$ and the winding number around $0$ of
the restriction of $\gamma$ to $[t_1,t_2]$ is
$1$.  The only reason that $\kappa=8$ is excluded from the theorem
is that it has not been proven yet that the SLE path $\gamma$ is
a continuous path when $\kappa=8$; this was proved
for all other $\kappa$ in~\cite{RS}.  In the range
$\kappa\in[0,4]$ the theorem holds trivially with $\lambda=0$, because
 $\gamma$ is a.s.\ a simple path.

\medskip\noindent{\bf Acknowlegments.}
This paper was planned prior to Smirnov's~\cite{Spaper}.
At the time, the link with discrete percolation was
only conjectural.
The conjecture has been established in~\cite{Spaper},
and therefore Theorems~\ref{t.pexpo} and \ref{t.contexpo} 
can now be stated as unconditional theorems.
Therefore, a significant portion of the credit for these results
also belongs to S.~Smirnov.

We wish to thank Harry Kesten and Michael Solomiak for
useful advice.

\section{The scaling limit exponent}
\begin{figure} 
\SetLabels
\B(.78*.96)$A_\theta$\\
\endSetLabels
\centerline{\AffixLabels{\includegraphics*[height=2.4in]
{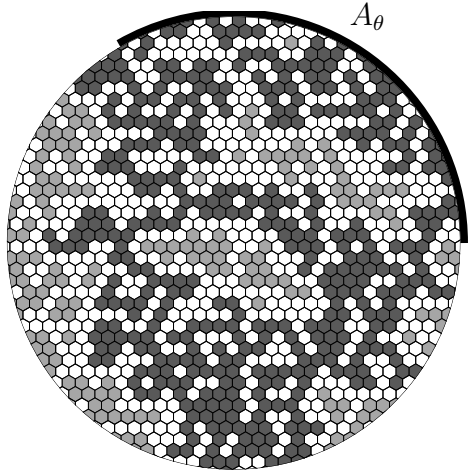}}}
\caption{\label{f.F}The set $Q_\delta(\theta)$.}
\end{figure}

Let $\theta\in[0,2\pi]$, and let $A_\theta$ be the arc
$$
A_\theta:= \bigl\{e^{is}:s\in[0,\theta]\bigr\}\subset\p\U\,.
$$
Fix some $\delta>0$, and consider site percolation
with parameter $p=1/2$ on the triangular grid of mesh $\delta$. 
It is  convenient
to represent a percolation configuration by coloring
the hexagonal faces of the dual grid, black for open,
white for closed.  Let $\mathfrak B_\delta$ denote the
union of all the black hexagons, and let $Q_\delta(\theta)$ denote
the union of $A_\theta$ with all the connected components
of $\mathfrak B_\delta\cap\closure\U$ which meet $A_\theta$.
See Figure~\ref{f.F}. 
Let $Q(\theta)$ denote a random subset of $\closure{\U}$ whose
law is the weak limit as $\delta\downarrow0$ of the law of
$Q_\delta(\theta)$.
(The law of $Q_\delta(\theta)$ can be thought of as a probability
measure
on the Hausdorff space of compact subsets of $\closure\U$.)
By~\cite{Spaper,Sprep}, the limit exists, and, moreover, it can be
explicitly described via \sle6/, as follows.

Recall from~\cite{Spaper} that the limit $\tilde\gamma$
as $\delta\downarrow0$ of the
outer boundary of $Q_\delta(\theta)$ consists of $A_\theta$
concatenated with the path of chordal \sle6/ joining the endpoints
of $A_\theta$.  By~\cite{Spaper,Sprep}, the scaling limit of percolation
is conformally invariant.  Hence, it is not too hard to verify that
$Q(\theta)$ is obtained by ``filling in'' each component
of $\U\setminus\tilde\gamma$ with an independent copy
of $Q(2\pi)$.  A more precise statement is the following.

\begin{theorem}[Smirnov]\label{t.Qt}
Let $\theta\in(0,2\pi)$,
let $\gamma_\ch$ denote the path of chordal \sle6/ from $1$ to
$e^{i\theta}$
in $\U$, and let $\tilde\gamma$ denote the curve obtained by
concatenating $A_\theta$, clockwise oriented, with $\gamma_\ch$.
Let $\mathcal W$ denote the collection of all
connected components $W$ of $\U\setminus\tilde\gamma$
such that $\tilde\gamma$ has winding number $-1$ about points
in $W$.  For each $W\in\mathcal W$, let $\psi_W:\U\to W$
be a conformal homeomorphism, and let $Q^W$ denote an independent
copy of $Q(2\pi)$.  Then the distribution
of $Q(\theta)$ is the same as the distribution of
$$
\bigcup_{W\in\mathcal W} \psi_W(Q^W)\,.
$$
(This union contains $\tilde\gamma$ a.s.)
\end{theorem}

Theorem~\ref{t.Qt} follows from the results and methods
of~\cite{Spaper}.
See~\cite{Sprep} for further details.

\medskip
We are interested in the distribution of the distance from $Q(\theta)$
to $0$.
However, it is more convenient to study the 
distribution of a very closely related quantity,
the {\em conformal radius}. 
 Let $U(\theta)$ denote the component of $0$ in $\U\setminus
Q(\theta)$. (It follows from 
Theorem~\ref{t.Qt} or from the Russo-Seymour-Welsh Theorem that $0\notin
Q(\theta)$ a.s.) Let $\psi=\psi_\theta:U(\theta)\to\U$ be the conformal
map,
normalized by $\psi(0)=0$ and $\psi'(0)>0$.
Define the conformal radius $\crad(\theta)$
of $U(\theta)$ about $0$ by $\crad(\theta):=1/\psi'(0)$.  
A well known consequence of the Koebe 1/4 Theorem and the Schwarz Lemma
(see, e.g.,~\cite{A}) is that
\begin{equation}\label{e.k1q}
\frac {\crad(\theta)}4\le\dist\bigl(0,Q(\theta)\bigr)\le
\crad(\theta)\,,
\end{equation}
and therefore information about the distribution of
$\crad(\theta)$ translates to information about the distribution
of $\dist\bigl(0,Q(\theta)\bigr)$.
Set
\begin{equation}\label{e.hdef}
h(\theta,t):= \Pb{ \crad(\theta)\le e^{-t}}\,.
\end{equation}
Note that for $t>0$,
\begin{equation}\label{e.easybd}
h(0,t)=\lim_{\theta\downarrow0} h(\theta,t)=0\,,
\end{equation}
follows  from the Russo-Seymour-Welsh Theorem
and the observation that
$\crad(\theta)$ tends to $1$ if the diameter of
$Q(\theta)$ tends to zero.

\begin{lemma}\label{l.pde}
In the range $\theta\in(0,2\pi)$, $t>0$,
the function $h$ satisfies the PDE
\begin{equation}\label{e.pde}
\frac\kappa2\, \p_\theta^2 h(\theta,t)
+\cot\bigl(\frac\theta2\bigr)\; \p_\theta h(\theta,t)
-\p_t h(\theta,t)=0
\,,
\end{equation}
with $\kappa=6$.
\end{lemma}

(The reason that we keep $\kappa$ as a parameter is that some
of the work done below will also apply to the proof of
Theorem~\ref{t.sletwist}, where $\kappa$ is not necessarily $6$.)

\proof
We assume that $\gamma_\ch$ and $Q(\theta)$ are coupled as
in Theorem~\ref{t.Qt}.
Let $T_\ch=T_\ch(\theta)$ be the first time $t$ such that
$\gamma_\ch[0,t]$ disconnects $0$ from $e^{i\theta}$ in $\closure\U$.
Let $\gamma_\ra$ denote the path of radial \sle6/
from $1$ to $0$ in $\closure\U$,
and let $T=T_\ra(\theta)$ denote the first time $t$ such that
$\gamma_\ra[0,t]$ disconnects $0$ from $e^{i\theta}$ in $\closure\U$.
{}From~\cite[Theorem 4.1]{LSW2} we know that up to time
reparameterization,
the restriction of $\gamma_\ra$ to $[0,T]$
has the same distribution as the restriction of
$\gamma_\ch$ to $[0,T_\ch]$.  
We therefore assume that indeed $\gamma_\ra$ restricted to $[0,T]$
is a reparameterization of $\gamma_\ch$ restricted to $[0,T_\ch]$.
Let $U_t$ be the connected component of $\U\setminus\gamma_\ra[0,t]$
which contains $0$, and let
$$
g_t:U_t\to\U
$$
be the conformal map, normalized by $g_t(0)=0$, $g_t'(0)>0$.
By the definition of radial \sle6/, the maps $g_t$ satisfy
\begin{equation}\label{e.radial}
\p_t g_t(z) = - g_t(z)\,
\frac {g_t(z)+e^{i\sqrt\kappa B_t}}
{g_t(z)-e^{i\sqrt\kappa B_t}}\,,\qquad g_0(z)=z\,,
\end{equation}
where $\kappa=6$ and $B_t$ is Brownian motion on $\R$ with
$B_0=0$.
Differentiating~\eref{e.radial} with respect to
$z$ at $z=0$ gives
\begin{equation}\label{e.derz}
g_t'(0)=e^{t}\,.
\end{equation}

At time $T$, $\gamma_\ra$ separates $0$ from
$e^{i\theta}$.
Let $W'$ be the connected component of
$\U\setminus\gamma_\ra[0,T]$ which contains $0$.
  We distinguish between two possibilities.
Either the boundary of $W'$ on $\gamma_\ra$ is
on the right hand side of $\gamma_\ra$, in which
case set $\nu:=1$, or on the left hand side,
and then set $\nu:=-1$.
In the notation of Theorem~\ref{t.Qt},
if $\nu=-1$, then $W'\notin\mathcal W$.
Hence, by~\eref{e.derz} and the coupling of $\gamma_\ra$ with
$Q(\theta)$,
\begin{equation}\label{e.cc}
\crad(\theta) = e^{-T}\,,\qquad\hbox{if }\nu=-1\,.
\end{equation}

We also want to understand the distribution of $\crad(\theta)$
given $\gamma_\ra$ and $\nu=1$.  In that case, $W'\in\mathcal W$.
In Theorem~\ref{t.Qt}, we may take the map $\psi_{W'}$
to satisfy $\psi_{W'}(0)=0$.
Let $r'$ denote the conformal radius about $0$ of
$\U\setminus Q^{W'}$, where $Q^{W'}$ is as in Theorem~\ref{t.Qt}.
Then $r'$ has the same distribution as $\crad(2\pi)$
and is independent from $\bigcup_{s>0}\mathcal F_s$,
where $\mathcal F_s$ denotes the $\sigma$-field generated
by $(B_t:t\le s)$. 
Moreover, by the chain rule for the derivative at zero,
\begin{equation}\label{e.ccm}
\crad(\theta) = r'\,e^{-T}\,,\qquad\hbox{if }\nu=1\,.
\end{equation}

To make use of the relations~\eref{e.cc} and~\eref{e.ccm}, we have
to understand the relation between $\nu$ and $\mathcal F_T$.
Recall that
$g_t\bigl(\gamma_\ra(t)\bigr)=\exp(i\sqrt\kappa B_t)$.
Set for $t<T$,
\begin{equation}\label{e.Ytdef}
Y_t=Y_t^{\theta}:=-i\log g_t(e^{i\theta})-\sqrt\kappa\, B_t\,,
\end{equation}
with $Y_0^\theta=\theta$ and the $\log$ branch chosen continuously.
That is, $Y_t$ is the length of the arc on 
$\p\U$ which corresponds under $g_t^{-1}$ to
the union of $A_\theta$ 
with the right hand side $\gamma_{\rh}[0,t]$ of $\gamma_\ra[0,t]$.
Suppose, for a moment that $\nu=1$.  That means that the boundary
of $W'$ is contained in $\gamma_{\rh}$
and $A_\theta$.  Consequently, as $t\uparrow T$, the harmonic
measure from $0$ of $A_\theta
\cup \gamma_\rh$ tends to $1$.  By conformal invariance
of harmonic measure,  this means that
$Y_t\to 2\pi$ as $t\uparrow T$ on the event $\nu=1$.
Similarly, we have
$Y_t\to 0$ as $t\uparrow T$ on the event $\nu=-1$.
Set $Y_T:=\lim_{t\uparrow T} Y_t$.
By~\eref{e.cc} and~\eref{e.ccm}, we now have
$$
\PB{\crad(\theta)\le e^{-t}\md \mathcal F_T}
=
1_{\{Y_T=2\pi\}}\,\Pb{r'\le  e^{T-t}\md \mathcal F_T}\,
+
1_{\{Y_T=0,\,t\le T\}}\,.
$$
Taking expectation gives
\begin{equation}
\label{e.heq}
h(\theta,t) = \Eb{h(Y^\theta_T,t-T)}\,,
\end{equation}
where we use the fact that $h(\theta, t) =1$ for $t\le 0$
and $\theta \in [0,2\pi]$.
{}From~\eref{e.radial} and~\eref{e.Ytdef} we get 
\begin{equation}\label{e.dY}
dY_t= \cot(Y_t/2)\, dt -\sqrt\kappa\,dB_t\,.
\end{equation}
By~\eref{e.heq}, for every constant $s>0$ the process
 $t\mapsto h(Y^\theta_t,s-t)$ is a local martingale on
$t<\min\{T,s\}$.
The theory of diffusion processes and~\eref{e.heq},~\eref{e.dY}
imply that $h(\theta,t)$ is smooth in the range $\theta\in(0,2\pi)$,
$t>0$.  
By It\^o's formula at time $t=0$
$$
dh(Y_t,s-t)=\Bigl(\frac \kappa 2\, \p_\theta^2 h(\theta,s)
+\cot\frac\theta2\; \p_\theta h(\theta,s)
-\p_s h(\theta,s)\Bigr)dt -\sqrt\kappa\,\p_\theta h(\theta,s)\, dB_t.
$$
Since $h(Y_t,s-t)$ is a local martingale, the $dt$ term must
vanish, and we obtain~\eref{e.pde}.
\QED

In order to derive boundary conditions, 
it turns out that it is more convenient to work with the following
modified version of $h$:
$$
\th(\theta,t):= \int_0^1 h(\theta,t+s)\,ds\,.
$$
Since, $h$ is smooth, $\tilde h$ also satisfies the PDE~\eref {e.pde}.
\begin{lemma}\label{l.bdval}
For every fixed $t>0$, the one sided $\theta$-derivative of
$\th$ at $2\pi$ is zero; that is,
\begin{equation}\label{e.bdval}
\lim_{\theta\uparrow 2\pi}
\frac{\th(2\pi,t)-\th(\theta,t)}{2\pi-\theta}=0\,.
\end{equation}
\end{lemma}

\proof
Let $\epsilon>0$ be very small
and set $\theta=2\pi-\epsilon$.
Let $\delta>0$ be smaller than $\epsilon$,
and let $\mathcal Z(r,\epsilon)$ denote the event that
there is a connected component
of $\mathfrak B_\delta\cap\closure\U$
which does not intersect $A_\theta$
but does intersect the two circles of radii $\epsilon$
and $r$ about the point $1$.

We claim that there are  constants $c,\alpha>0$ such that
\begin{equation}\label{e.Fr}
\Pb{\mathcal Z(r,\epsilon)} \le c (\epsilon/r)^{1+\alpha}\,,
\end{equation}
provided $0<\delta<\epsilon<r<2$.
This well known result is 
also used in Smirnov's arguments.
As we could not track down an explicit proof of this 
statement in  the literature and for the reader's convenience, 
we include a proof of this fact in the appendix.

Let $Q':= Q(2\pi)\setminus Q(\theta)$.  By letting $\delta\downarrow0$,
it follows from~\eref{e.Fr} that
\begin{equation}\label{e.diamest}
\Pb{\diam Q' \ge r} \le c\,(\epsilon/r)^{1+\alpha}\,.
\end{equation}
It is easy to verify that
there is a constant $c_1>0$ such that for every connected compact
$K\subset \closure\U$ which intersects $\p\U$, the harmonic measure
$\mu(\U,K,0)$ in $\U$ of $K$ from $0$ satisfies
\begin{equation}\label{e.diamharm}
c_1^{-1}\,\diam K\le \mu(\U,K,0)\le c_1\,\diam K\,.
\end{equation}
Moreover, if $\crad(K)$ denotes the conformal radius
of $\U\setminus K$, then
\begin{equation}\label{e.crad}
\min\Bigl\{-\log\crad(K),1\Bigr\}\le c_2 \bigl(\diam K\bigr)^2\,,
\end{equation}
where $c_2$ is some constant.
To justify this, note that $\crad(K)$ is monotone non-increasing in $K$
and hence $\crad(K)\ge \crad(B)$, where $B$ is the
smallest disk with $\p B$ orthogonal to $\p\U$ which contains
$K$.  For such a $B$, one can calculate $\crad(B)$ explicitly,
since the normalized conformal map from $\U\setminus B$ to
$\U$ is conjugate to the map $z\to \sqrt z$ by M\"obius transformations.
Now
\begin{eqnarray*}
\min\Bigl\{\log \crad(\theta)-\log\crad(2\pi),1\Bigr\}
&=&
\min\Bigl\{-\log\crad \bigl(\psi^{-1}_\theta(Q')\bigr),1\Bigr\}
\\
&\le &c_2 \bigl(\diam \psi^{-1}_\theta(Q')\bigr)^2
\\
&\le &c_2\,c_1^2\, \mu(\U,\psi^{-1}_\theta(Q'),0)^2
\\
&=& c_2\,c_1^2\, \mu(\U\setminus Q(\theta),Q',0)^2
\\
&\le & c_2\,c_1^2\, \mu(\U,Q',0)^2
\le c_2\,c_1^4\, (\diam Q')^2\,.
\end{eqnarray*}
(The equality is due to conformal invariance of harmonic
measure.)
Combining this with~\eref{e.diamest} gives
$$
\Pb{
\log \crad(\theta)-\log\crad(2\pi)
\ge s}
\le \min \bigl\{ 1, c_3\,(\epsilon/\sqrt s)^{1+\alpha}\bigr\} \,,\qquad
s\le 1\,.
$$
Hence, 
\begin{equation*}
\th(2\pi,t)-\th(\theta,t)
\le
\int_0^1
\Pb{
\log \crad(\theta)-\log\crad(2\pi)
\ge s}\,ds = o(\epsilon)\,,
\end{equation*}
and the lemma is established.
\QED

\proofof{Theorem \ref{t.contexpo}}
We are going to give a non-probabilistic proof 
based on the PDE and boundary conditions that we derived (but other 
justifications are also possible).
Let $\Lambda$ denote the differential operator on the left hand
side of~\eref{e.pde}, and set $S=(0,2\pi)\times(0,\infty)$.
Set
$$
q:= \frac{\kappa - 4}{\kappa}
\,,\qquad
H(\theta,t):=
\bigl(\sin(\theta/4)\bigr)^q\,  e^{-\lambda t}\,,
$$
where $\lambda$ is defined by~\eref{e.lqdef}.
The function $\th$ satisfies $\Lambda\th=0$ in $S$,
the Neumann boundary condition~\eref{e.bdval} the Dirichlet
condition $\th(0,t)=0$ for $t>0$, $0\le\th(\theta,t)\le1$
and $\th(\theta,0)>0$ for $\theta\in(0,2\pi]$.
Note that $H$ also has all these properties.
We now use the Maximum Principle to show that
there are positive constants $c, c'$ such that
\begin{equation}\label{e.Hsand}
\forall t\ge 1\,,\ \forall\theta\in[0,2\pi]\,,
\qquad
c\,  H(\theta,t)\le\th(\theta,t)\le c'\, H(\theta,t)
\,.
\end{equation}
This is a bit tricky because of the singularities of~\eref{e.pde}
at the boundary, but otherwise quite standard.   

Suppose that a function $G$ is defined and
continuous on $S^* :=\overline S \setminus \{ (0,0)\}$,
non-negative in some neighborhood of $(0,0)$ in $S^*$,
$G \ge 0$ on $t=0$, $G=0$ on $\theta=0$,
$\partial_\theta G = 0$ on $\theta=2 \pi$ and
$\Lambda (G) = 0$ in $S$.
Then, we claim that $G \ge 0$ on $S^*$.

Let $\epsilon>0$ and set
$
F:=G +\epsilon\,\theta^{q/2}+\epsilon\,.
$
Assume that $F<0$ somewhere in $S$.
The function $F$ is continuous on $S^*$, positive
in a neighborhood of $(0,0)$ in $S^*$ and
$F \ge \epsilon$ when $t=0$ or $\theta=0$. 
Continuity then implies that there is a pair $(\theta_0,t_0)\in S^*$
such that $F(\theta_0,t_0)\le 0$ and $t_0$ is minimal among all
such points in $S^*$.
Clearly, $t_0>0$ and $\theta_0>0$.  Observe also that $\theta_0\ne 2\pi$,
because $\p_\theta F>0$ at $\theta=2\pi$.  Hence,
$(\theta_0,t_0)\in S$.
A computation using the explicit value of $q$
shows that $\Lambda(F)= \epsilon \Lambda ( \theta^{q/2} ) <0$ in $S$.
By the minimality of $t_0$, we have
$$
\p_\theta F(\theta_0,t_0)=0\,,
\qquad
\p_\theta^2 F(\theta_0,t_0)\ge 0\,,
\qquad
\p_t F(\theta_0,t_0)\le 0\,.
$$
But the definition of $\Lambda$ shows
that this contradicts $\Lambda F<0$.
Therefore $F\ge 0$ on $S^*$.
Since $\epsilon>0$ was arbitrary, this proves that $G \ge 0$.

Now, note that $\Lambda (2-2t-\theta^2)<0$ for
$\theta\in(0,\pi)$.
Let
$$
F_1:=
c_1\,H-\th+2-2t-\theta^2\,,\qquad
F_2 :=
c_2 \, \th - H + 2 - 2t - \theta^2\,,
$$
where the constants $c_1, c_2$ are chosen so that
$F_1>0$ and $F_2>0$  on $\{\pi\}\times[0,1]$
and on $[0,\pi]\times\{0\}$.  The above argument
applied to the functions $F_j$ with $S$ replaced by
$(0,\pi)\times(0,1)$, shows that
$c_1\,H(\theta,1)- \th(\theta,1) \ge F_1 (\theta, 1) \ge 0$
and $c_2 \, \th (\theta, 1) - H (\theta, 1) \ge F_2 (\theta, 1) \ge 0$
for $\theta\in[0,\pi]$.   
By changing the constants $c_1$ and $c_2$ if necessary,
we may make sure that this same inequalities
hold for all $\theta\in[0,2\pi]$.
Yet another application of the same Maximum Principle argument,
this time in the range $[0, 2 \pi] \times [1, \infty)$ proves
that $c_1 \, H - \th \ge 0$ and $c_2 \, \th - H \ge 0$ for
$t\ge 1$, thereby establishing~\eref{e.Hsand}.  

The theorem now follows from~\eref{e.Hsand},
the definition of $\th$ and~\eref{e.k1q}.
\QED

\section{The discrete exponent}

We now show that the discrete exponent and the continuous
exponent we derived are the same.  The proof is quite direct.
However, for exponents related to crossings of several colors
(that is, black and white or occupied and vacant),
the analogous equivalence is not as easy.  
See~\cite{SmW} for a treatment of this more involved situation.

\proofof{Theorem \ref{t.pexpo}}
Let $u(r)$ denote the probability that
$Q=Q(2\pi)$ intersects the circle $r\p\U$,
and let $\mathcal C(R_1,R_2)$ denote the event 
that in site percolation with parameter $p=1/2$
on the standard triangular grid, there is an open cluster
crossing the annulus with radii $R_1$ and $R_2$ about $0$.
Let $\mathcal A(R)$ denote the event that there
is an open path separating the
circles of radii $R$ and $2R$ about $0$.
By the definition of the scaling limit it follows
that for all fixed $r\in(0,1)$ 
there is some $s_0=s_0(r)$ such that
\begin{equation}\label{e.compa}
\forall R\ge s_0,\qquad u(r/2)/2
\le
\Pb{\mathcal C(rR,R)}\le
2\,u(2r)\,.
\end{equation}
Let $R$ be large, $r>0$ small, set
$\tilde r = 2r$,
and let $N$ be the minimal integer satisfying $s_0(r)\,{\tilde
r}^{-N}>R$.
Note that the event $\{0\leftrightarrow \CS_R\}$ that
$0$ is connected to the circle $\CS_R$ satisfies
$$
\{0\leftrightarrow \CS_{R} \}\supset
\{0\leftrightarrow \CS_{2s_0}\}\cap \bigcap_{j=0}^{N-1}
\bigl(\mathcal C(s_0\, {\tilde r}^{-j}, 2\,s_0\,{\tilde r}^{-j-1})
\cap \mathcal A(s_0\,{\tilde r}^{-j})\bigr)\,.
$$
By the Harris/Fortuin-Kasteleyn-Ginibre (FKG) 
inequality, we therefore have
$$
\P[0\leftrightarrow \CS_R]
\ge \P[0\leftrightarrow \CS_{2s_0} ]
\prod_{j=0}^{N-1}
\Pb{ \mathcal C(s_0\, {\tilde r}^{-j}, 2\,s_0\,{\tilde r}^{-j-1})}
\prod_{j=0}^{N-1}
\Pb{\mathcal A(s_0\,{\tilde r}^{-j})}\,.
$$
By the 
Russo-Seymour-Welsh Theorem (RSW), there is a constant $c_1>0$ such
that $\Pb{\mathcal A(R)}\ge c_1$ for every $R$.  Applying this
and~\eref{e.compa} in the above, gives
$$
\P[0\leftrightarrow \CS_R]
\ge  \P[0\leftrightarrow \CS_{2s_0}]\,
c_1^N\, \bigl(u(r/2)/2\bigr)^N\,.
$$
By Theorem~\ref{t.contexpo}, $\log u(r)/\log r\to 5/48$
as $r\downarrow 0$.
Therefore,
$$
\liminf_{R\to\infty}\frac{\log \P[0\leftrightarrow \CS_R]}{\log R}
\ge\frac{\log (c_1/2)+\log u(r/2)}{-\log r}\to -5/48\,,
$$
as $r\downarrow 0$.  This proves the required lower bound
on $\P[0\leftrightarrow \CS_R]$.  The proof for the upper bound is
similar
(using independence on disjoint sets, in place of FKG),
 and is left to the reader.
\QED

\section{Counterclockwise loops for \sle\kappa/}

\proofof{Theorem~\ref{t.sletwist}}
Let $\kappa>4$, $\kappa\ne 8$.
Let $\gamma$ denote the path of radial \sle\kappa/
started from $1$, and for $\theta\in[0,2\pi]$ let 
$\gamma_\theta^t$ denote the concatenation of the counterclockwise
oriented arc $\p\U\setminus A_\theta$ with $\gamma[0,t]$.
Let $\mathcal E(\theta,t)$ denote the event
that $\gamma_\theta^t$ does not contain
a counterclockwise loop around $0$,
and set
$$
\hat h(\theta,t):=\Pb{\mathcal E(\theta,t)}.
$$
We know that a.s.\ for all
$t>0$, $\theta>0$ the path $\gamma[0,t]$ intersects
$A_\theta$.  (This follows by comparing
with a Bessel process. See, e.g.,~\cite{RS}.)
Hence, $\lim_{\theta\downarrow 0}\hat h(\theta,t)=0$
for all $t>0$.

Let $K_t$ denote the SLE hull
(i.e., $\gamma[0,t]$ union with the set
of points in $\closure\U$ separated
from $0$ by $\gamma[0,t]$), and
let $U_t:=\U\setminus K_t$,
which is the component of
$\U\setminus\gamma[0,t]$ containing $0$.
Let 
$$
I_t:=\bigl\{\theta\in(0,2\pi):e^{i\theta}\notin K_t\bigr\}\,.
$$
Let $g_t$ denote the conformal
map $g_t:U_t\to\U$ normalized by
$g_t(0)=0$ and $g_t'(0)>0$.

Let $T(\theta)$ denote the first time such
that $e^{i\theta}\in K_t$,
and let $\hat T:=\sup_\theta T(\theta)$.
Then $\hat T$ is also the first time $t$ such that
$\gamma[0,t]$ contains a loop around $0$,
and is the first time such that $I_t=\emptyset$.
Define $Y_t^\theta$ for $t<T(\theta)$
as in~\eref{e.Ytdef}.
A moment's thought shows that~\eref{e.heq}
holds with $\hat h$ in place of $h$.
Therefore, $\hat h$ satisfies~\eref{e.pde}
for $\theta\in(0,2\pi)$ and $t>0$.

We now need to verify that $\hat h$ also
satisfies the Neumann boundary condition~\eref{e.bdval}.
It is immediate that there is a constant
$c_1>0$ such that 
\begin{equation}\label{e.Harn}
\forall t>0,\ \forall
\theta\in[\pi/2,3\pi/2],\qquad
\bigl|\p_\theta \hat h(\theta,t)\bigr|\le c_1
\,.
\end{equation}
(This is a Harnack type inequality.)
For $\theta\in[0,2\pi]$ and $t\in[T(\theta),\hat T)$ {\em define},
\begin{equation*}
      Y_t^\theta :=
\begin{cases} \inf\{Y_t^\alpha: \alpha\in I_t \}& \theta\le\inf I_t\,,
\\
\sup\{Y_t^\alpha: \alpha\in I_t \}& \theta\ge\sup I_t\,.
\end{cases}
\end{equation*}
For fixed $\theta$, the process $Y_t^\theta$ is a process
satisfying~\eref{e.dY} with instantaneous reflection
in the boundary (see, e.g., \cite{RY} for a detailed
treatment of such processes) stopped at the time $\hat T$.
Clearly, for all $t\in[0,\hat T)$ the set of $\theta\in[0,2\pi]$ for
which
$\mathcal E(\theta,t)$ holds is 
$\bigl\{\theta\in[0,2\pi]:\inf_{s\le t} Y^\theta_s>0\bigr\}$.
Also, for all $t\ge \hat T$ the set of $\theta\in[0,2\pi]$
for which $\mathcal E(\theta,t)$ holds is either $\emptyset$
or 
$\bigl\{\theta\in[0,2\pi]:\inf_{s< \hat T} Y^\theta_s>0\bigr\}$.

Now let $\theta_1$ and $\theta_2$ satisfy
$3\pi/2<\theta_1<\theta_2<2\pi$,
and fix some $t_1>0$.
Let $\tau:=\inf\bigl\{t\in [0,\hat T]: Y^{\theta_1}_t=\pi\bigr\}$,
and let $\mathcal S$ be the event that $\tau\le t_1$ and
$Y^{\theta_2}_\tau\ne\pi$.
Observe that $\mathcal E(\theta_1,t_1)$
is equivalent to $\mathcal E(\theta_2,t_1)$ on the complement of
$\mathcal S$.  Consequently, the strong Markov property at time $\tau$
gives
\begin{equation}\label{e.hdif}
\hat h(\theta_2,t_1)-\hat h(\theta_1,t_1)
\le \P[\mathcal S]\,
\Eb{\hat h(Y^{\theta_2}_\tau,t_1-\tau)-
\hat h(Y^{\theta_1}_\tau,t_1-\tau)\md \mathcal S}\,.
\end{equation}
Observe from~\eref{e.dY} that
$$
\p_t\bigl( Y^{\theta_2}_t-Y^{\theta_1}_t\bigr)\le 0
\,.
$$
Therefore, on $\mathcal S$ we have
$Y^{\theta_2}_\tau- Y^{\theta_1}_\tau\le \theta_2-\theta_1$.
Hence, by~\eref{e.Harn} combined with~\eref{e.hdif},
\begin{equation}\label{e.hd}
\hat h(\theta_2,t_1)-\hat h(\theta_1,t_1)
\le c_1\, \P[\mathcal S]\, (\theta_2-\theta_1)\,.
\end{equation}
But when $\theta_1<\theta_2\le 2\pi$ are both close to $2\pi$,
we know that with probability close to $1$ there will
be some time $t_0\in(0,t_1)$ such that $Y^{\theta_1}_{t_0}=2\pi$.
In that case also $Y^{\theta_2}_{t_0}=2\pi$, and hence
$Y^{\theta_1}_{t}=Y^{\theta_2}_{t}$ for all $t\in[ t_0,\hat T)$.
In particular, $\P[\mathcal S]\to 0$ as $\theta_1\to2\pi$. 
By~\eref{e.hd} we get the Neumann condition:
$$
\lim_{\theta\uparrow2\pi}\p_\theta\hat h(\theta,t_1)=0\,.
$$
Now the proof of the theorem is completed exactly as for
the corresponding proof of~\eref{e.Hsand}.
\QED

\appendix\section{Appendix: an a priori half plane exponent}
We now give a brief outline of the proof of the a priori estimate (\ref
{e.Fr})
on the probability that there exists three disjoint crossings of a 
half-annulus.  This is based on the exponent for the existence of
two disjoint crossings.

\begin{lemma}\label{l.twocross}
There is a constant $c>0$ such that for all $R>r>2$ the probability
$f(R,r)$ that there are two disjoint crossings between
$\CS_R$ and $\CS_r$ within the half-plane
$i\closure\H=\bigl\{z:\Re z\le0\}$ in site percolation with
parameter $p=1/2$ on the standard triangular grid
satisfies
$$
c^{-1} \,r/R\le f(R,r)\le c\,r/R\,.
$$
\end{lemma}

Loosely speaking, this lemma says that the two-arms
half plane exponent is $1$. 
Note that the half-plane exponents for \sle6/ were calculated
in~\cite{LSW1,LSW3}.  They are now valid for critical site percolation
on the triangular grid, since \sle6/ describes the scaling limit.
However, to show the equivalence, this lemma seems to be needed.

\proof
Let $S_R$ be the strip $S_R:=\bigl\{z\in\C:\Re z\in [-R,0]\bigr\}$,
and let $X_R$ be the set of clusters in $S_R$ which join
the interval $I_R:= [-iR,iR]\subset i\R$
to the line $L_R:= -R + i\R $.  It follows from RSW that
\begin{equation}\label{e.XR}
\sup_{R>0}\Eb{|X_R|}<\infty\,,
\qquad
\inf_{R>2}\Pb{|X_R|>2}>0\,.
\end{equation}
Define the event $\mathcal D_R$ that 
there exists a path of occupied site in $S_R$ joining
the origin to $L_R$ and a path of vacant sites in $S_R$ joining
the vertex directly below the origin to $L_R$.
Note that $\mathcal D_R$ is the same as the event
that there is some $C\in X_R$ such that
$0$ is the lowest vertex in $C\cap i\R$.
By invariance under vertical translations,
since the number of vertices in $I_R$
is of order $R$,~\eref{e.XR} implies that
$$
\forall R>2,\qquad
\Prob [\mathcal D_R ]\asymp 1/R\,,
$$
where $g\asymp h$ means that $g/h$ is bounded and
bounded away from zero.

We now employ the standard ``color-flipping'' argument
(see e.g., \cite {ADA,LSW1}).
That is, suppose that the origin is connected to $L_R$
in $S_R$.  By flipping all the vertices
below the topmost crossing from $0$ to $L_R$ in $S_R$,
it follows that $\P[\mathcal D_R]$ is 
equal to the probability of the event
$\mathcal D^2_R$ that there are disjoint open
paths in $S_R$ connecting the origin and the vertex
below the origin to $L_R$.
Let $\mathcal C^2_R$ be the event that there are
two disjoint paths connecting the origin and
the vertex below the origin to $\CS_R$ in $i\closure\H$.
Clearly, $\mathcal C^2_R\supset \mathcal D^2_R$.
On the other hand, when $R$ is sufficiently large,
by RSW the probability that
there are two disjoint left-right crossings of
each of the the rectangles $[-R,0]\times [R/10,R/5]$
and $[-R,0]\times [-R/5,-R/10]$ is bounded away
from zero.
Therefore,
the FKG inequality implies that there is some constant $c>0$ such that
$$
\P[\mathcal D^2_R]\ge c\, \P[\mathcal C^2_{2R}]\,.
$$
Consequently, we find that $\P[\mathcal C^2_R]\asymp R^{-1}$ for $R>2$.
Another application of RSW and FKG
as in the proof of Theorem~\ref{t.pexpo} shows that
$$
\P[\mathcal C^2_R]\asymp f(R,r)\,\P[\mathcal C^2_r]\,,
\qquad\hbox{when } R>2r>4\,.
$$
The lemma follows. 
\QED

By using RSW yet again combined with the van den Berg-Kesten inequality
(BK),
it follows that there are $c,\epsilon>0$
such that the probability that $\CS_r$ and $\CS_R$ are
connected in $i\closure\H$ by two open crossings and one closed
crossing is at most $c\,(r/R)^{1+\epsilon}$, when $r>1$.
This immediately implies~\eref {e.Fr}, because the
existence of a cluster in $\closure\U$ intersecting
$\CS_r$ and $\CS_\epsilon$ but not intersecting
$A_\theta$ implies that there are also two disjoint closed crossings
from $\CS_\epsilon$ to $\CS_r$ in $\closure\U$.
\bigskip

By looking at the right-most vertices of clusters of diameter greater than
$R$,
an argument similar to the proof of Lemma~\ref{l.twocross}
can be used to prove directly that the three-arms exponent
in a half-plane is $2$.  
An analogous argument \cite[Lemma 5] {KSZ}
shows that the multicolor five-arms exponent 
in the plane is equal to $2$, which then implies
that the multicolor six arms  exponent is strictly larger than
$2$. (This is also an instrumental a priori
estimate in~\cite {Sprep}.)
Similarly, the Brownian intersection exponent $\xi (2,1)=2$ can
be easily determined~\cite {Lbook}, and a
direct proof also works for some exponents associated 
to loop-erased random walks (e.g.,~\cite{Lkes}).
However, these exponents that take the value $1$ or $2$ are 
exceptional. 
For the determination
of  fractional exponents, such as, for instance, the 
one-arm exponent in the present paper,  
such direct arguments do not work. 

\section{The monochromatic 2-arm exponent}\label{a.twoarm}

This informal appendix will discuss the monochromatic 2-arm 
exponent (sometimes also called ``backbone exponent'');
 that is, the exponent describing the decay
of the probability that there are two disjoint open
crossings from $\CS_r$ to $\CS_R$.  So far, despite some
attempts, there is no established prediction in the theoretical 
physics community for  
the value of
this exponent, as far as we are aware.
Simulations \cite {Gr} show that its numerical value is 
$.3568 \pm .0008$.

The color switching argument used above shows that for
half-plane exponents, the required colors of the arms do not
matter.  For full plane exponents, this argument can
only show that every two $k$-arm exponents are
the same, if both colors are present in both \cite{ADA}.
For example, the plane exponent for the existence
of open, closed, open, closed, open crossings,
in clockwise order, is the same as the exponent for
the closed, open, open, open, open exponent.
Thus, for each $k=2,3,\dots$, there are two types
of full-plane exponents, the monochromatic $k$-arm exponent,
and the multichromatic $k$-arm exponent.  The multichromatic
exponents with $k\ge 2$ can all be worked out 
from~\cite{LSW2} combined
with~\cite{Spaper} (see \cite {SmW}).
We now show  (ommiting many details)
that the monochromatic $2$-arm exponent is in
fact also the leading eigenvalue 
of a differential operator.  
So
far, we
have not been able to solve explicitely the PDE problem and to give
an explicit formula for the value of the exponent, but 
this might perhaps be doable by someone more proficient in this field.

Like in the one-arm case, we have to work with a slightly
more general problem. 
Let $\alpha,\beta>0$ satisfy $\alpha+\beta<2\pi$,
and let $A$ and $A'$ be two consecutive arcs on
$\p\U$ with lengths $\alpha$ and $\beta$, respectively, and
set $A''=\p\U\setminus (A\cup A')$.
Let $\gamma = 2 \pi -  \alpha - \beta$ denote the length of $A''$
and let $Q$ be the scaling limit of the set of points $p\in\closure\U$
which are connected to $A''$ by two paths lying in the
union of $A'$ with the set of black hexagons,
such that the two paths do not share any hexagon
except for the hexagon containing $p$.
Let $\kk(\alpha,\gamma,t)$ be the probability that the conformal
radius of $\U\setminus Q$ about $0$ is at most $e^{-t}$.
Using arguments as above,
it is easy to check, by starting an \sle6/ path from the common
endpoint of $A$ and $A'$, that $\kk(\alpha,\gamma,t)$
satisfies the PDE
$$
3 \,  \p_\alpha^2 \kk  
+ \cot\bigl(\frac \alpha 2\bigr) \p_\alpha \kk
+\Bigl(\cot\bigl(\frac {\alpha + \gamma} 2\bigr)
- \cot\bigl(\frac \alpha 2 \bigr) \Bigr) \p_\gamma \kk - \p_t \kk = 0.
$$
The Dirichlet condition
$$
\kk=0\,,\qquad\hbox{when }
\gamma = 0$$
is immediate. 
Similarly, one can check that for all $\gamma \not= 0$, 
$$ 
\kk ( 0, \gamma ) = \kk ( 0 , 2 \pi ).
$$
Finally, using arguments as for the one-arm
exponent, one can show that 
the boundary condition on $\alpha + \gamma = 2 \pi$ is 
$$
\partial_\gamma \kk = 0 .
$$
Hence, the monochromatic
two-arm exponent is a number (and in fact the  
 unique number)  
$\lambda>0$ such
that there is a non-negative function $\kk_1(\alpha,\gamma)$,
not identically zero, such
that $e^{-\lambda t}\, \kk_1(\alpha,\gamma)$ solves the above PDE and
has the corresponding boundary behaviour (one can fix 
the multiplicative constant by setting $\kk_1 ( 0, 2 \pi ) = 1$).
We have been unable to guess an explicit solution 
to this boundary value problem.
Some insight might be gained by noting that when 
$\alpha$ and $\gamma$ go down to zero, the solution should behave like
Cardy's~\cite{CaFormula} formula in $\alpha  / ( \alpha + \gamma )$,
which is a hypergeometric function.

One can rewrite the PDE using other variables. In terms of the 
function $ \kk_2( \alpha, \beta) := \kk_1 ( \alpha, 2 \pi - \alpha - \beta)$,
one looks for the non-negative eigenfunctions of the more symmetric operator 
$$
3 (\partial_\alpha^2 - 2 \partial_\alpha \partial_\beta 
+ \partial_\beta^2 )  
+ \cot\bigl(\frac \alpha 2\bigr) \p_\alpha 
- \cot\bigl(\frac \beta 2\bigr) \p_\beta 
$$
 with boundary conditions (for $\alpha , \beta>0$, $\alpha+\beta<2\pi$),  
$$
\kk_2 (\alpha, 2 \pi - \alpha) = 0 , \  
\kk_2 (0, \beta) = 1  
\hbox { and } 
\partial_\beta \kk_2 ( \alpha, 0) = 0 .$$

\bigskip

\filbreak
\begingroup
\small
\parindent=0pt

\vtop{
\hsize=2.3in
Greg Lawler\\
Department of Mathematics\\
Box 90320\\
Duke University\\
Durham NC 27708-0320, USA\\
{jose@math.duke.edu}
}
\bigskip
\vtop{
\hsize=2.3in
Oded Schramm\\
Microsoft Corporation,\\
One Microsoft Way,\\
Redmond, WA 98052, USA\\
{schramm@microsoft.com}
}
\bigskip
\vtop{
\hsize=2.3in
Wendelin Werner\\
D\'epartement de Math\'ematiques\\
B\^at. 425\\
Universit\'e Paris-Sud\\
91405 ORSAY cedex, France\\
{wendelin.werner@math.u-psud.fr}
}
\endgroup

\filbreak

\end {document}